\documentclass[preprint,12pt, a4paper]{elsarticle}

\usepackage{amssymb}
\usepackage{amsfonts}
\usepackage{lineno}
\usepackage{float}
\restylefloat{table}

\usepackage[T1]{fontenc}
\usepackage[utf8]{inputenc}
\usepackage{url}
\usepackage{tikz}
\usepackage[english]{babel}
\usepackage{xcolor}
\definecolor{pyblue}{rgb}{0,0,0.5}
\definecolor{pyred}{rgb}{0.6,0,0}
\definecolor{pygreen}{rgb}{0,0.5,0}
\definecolor{pygray}{rgb}{0.5, 0.5, 0.5}
\definecolor{lightyellow}{RGB}{255,253,233}
\usepackage{listings}
\usepackage{braket}
\usepackage[defaultlines=2,all]{nowidow}
\usepackage{booktabs}
\usepackage{subcaption}
\usepackage[retain-zero-exponent, per-mode=symbol, output-exponent-marker=e, group-digits=false]{siunitx}
\usepackage{nicefrac}
\usepackage{amsmath}
\usepackage{hyperref}

\newcommand{\transposed}{^\top}

\newcommand{\R}{\mathbb{R}}

\newcommand{\dmeas}[1]{\ \mathrm{d}#1}
\newcommand{\dual}[3]{\left\langle #1 , #2 \right\rangle_{#3}}

\newcommand{\grad}{\nabla}

\newcommand{\integral}[1]{\int_{#1}}

\newcommand{\admissiblegeom}{\mathcal{A}}

\newcommand{\costfunction}{\mathcal{J}}

\newcommand{\vectorfield}{\mathcal{V}}

\newcommand{\qe}[1]{``#1''}

\usetikzlibrary{positioning, arrows, backgrounds, shapes.multipart, calc}
\newlength\figureheight
\newlength\figurewidth

\DeclareFixedFont{\ttb}{T1}{txtt}{bx}{n}{10} 
\DeclareFixedFont{\ttm}{T1}{txtt}{m}{n}{10}  

\newcommand\pythonstyle{\lstset{
		language=Python,
		basicstyle=\ttm,
		otherkeywords={self},             
		keywordstyle=\ttb\color{pyblue},
		emph={MyClass,__init__},          
		emphstyle=\ttb\color{pyred},    
		stringstyle=\color{pygreen},
		commentstyle=\color{pygray},
		frame=single,                         
		showstringspaces=false,
		tabsize=4,
		numbers=left,
		backgroundcolor=\ttb\color{lightyellow}
}}

\lstnewenvironment{python}[1][]
{
	\pythonstyle
	\lstset{#1}
}
{}

\newcommand\pythoninline[1]{{\pythonstyle\lstinline!#1!}}

\journal{arXiv}

\begin{document}

\begin{frontmatter}

\title{cashocs: A Computational, Adjoint-Based Shape Optimization and Optimal Control Software}

\author[ITWM,TUK]{Sebastian Blauth}
\address[ITWM]{Fraunhofer ITWM, Kaiserslautern, Germany}
\address[TUK]{TU Kaiserslautern, Kaiserslautern, Germany}

\ead{sebastian.blauth@itwm.fraunhofer.de}

\begin{abstract}
	The solution of optimization problems constrained by partial differential equations (PDEs) plays an important role in many areas of science and industry. In this work we present cashocs, a new software package written in Python, which automatically solves such problems in the context of optimal control and shape optimization. The software cashocs implements a discretization of the continuous adjoint approach, which derives the necessary adjoint systems and (shape) derivatives in an automated fashion. As cashocs is based on the finite element software FEniCS, it inherits its simple, high-level user interface. This makes it straightforward to define and solve PDE constrained optimization problems with our software. In this paper, we discuss the design and functionalities of cashocs and also demonstrate its straightforward usability and applicability.
\end{abstract}


\end{frontmatter}

\section*{Required Metadata}
\label{req_meta}

\section*{Current code version}
\label{current_code}

\begin{table}[H]
\begin{tabular}{|l|p{6.5cm}|p{6.5cm}|}
\hline
\textbf{Nr.} & \textbf{Code metadata description} & \\
\hline
C1 & Current code version & v1.0.3 \\
\hline
C2 & Permanent link to code/repository used for this code version & \url{https://github.com/sblauth/cashocs/releases/tag/v1.0.3} \\
\hline
C3 & Code Ocean compute capsule & NA\\
\hline
C4 & Legal Code License   & GNU GPL v3.0 (or later) \\
\hline
C5 & Code versioning system used & git \\
\hline
C6 & Software code languages, tools, and services used & Python, FEniCS, NumPy, PETSc, meshio, Gmsh \\
\hline
C7 & Compilation requirements, operating environments \& dependencies & FEniCS, meshio, Gmsh, matplotlib \\
\hline
C8 & If available Link to developer documentation/manual & \url{https://cashocs.readthedocs.io/}\\
\hline
C9 & Support email for questions & \href{mailto:sebastian.blauth@itwm.fraunhofer.de}{sebastian.blauth@itwm.fraunhofer.de} \\
\hline
\end{tabular}
\caption{Code metadata (mandatory)}
\label{metadata} 
\end{table}

%

\section{Motivation and significance}
\label{sec:motivation}

Shape optimization and optimal control problems constrained by partial differential equations (PDEs) and their numerical solution are important in many areas of science and industry: They are, for example, used for the optimization of chemical reactors \cite{Blauth2020Optimal}, glass cooling processes \cite{Pinnau2004Optimal}, and semiconductors \cite{Hinze2002optimal} as well as the optimal design of cooling systems \cite{Blauth2020Shape}, aircrafts \cite{Schmidt2013Three}, and electric machines \cite{Gangl2015Shape}. To solve these problems, the so-called adjoint approach is often employed, which facilitates the computation of (shape) gradients for the problems, which can be used to solve them numerically. 
However, for complex, coupled, or highly nonlinear problems, such as the ones arising from industrial applications, even the derivation of the necessary equations for the adjoint approach is an extremely involved and error-prone task. Consequently, it is not feasible to carry out the adjoint approach manually anymore (see, e.g., \cite{Mitusch2019dolfin}). For these reasons, there has been a lot of effort recently to automate the tasks for solving PDE constrained optimization problems, resulting in software such as dolfin-adjoint \cite{Dokken2020Automatic} and Fireshape \cite{Paganini2020Fireshape}, shape optimization capabilities for the finite element software NGSolve \cite{Gangl2020Fully}, and our software cashocs. 

What distinguishes cashocs from these other packages is its novel approach of using automatic differentiation solely to derive the adjoint system and (shape) derivatives, while implementing and automating a discretization of the continuous adjoint approach in all remaining aspects. 
This means that the optimization algorithms together with all required operations are implemented as discretizations of the underlying infinite-dimensional operations. The aforementioned operations include, e.g., the determination of (shape) gradients from the computed (shape) derivatives, the discretization and numerical solution of the state and adjoint equations, the computation of scalar products, and the usage of projection operators.
Therefore, the calculated (shape) derivatives are only used as \qe{inputs} for our framework. Particularly, cashocs implements discretizations of continuous, infinite-dimensional optimization algorithms which are strongly related to the underlying optimization problem, whereas the other packages use either external optimization libraries \cite{Dokken2020Automatic, Paganini2020Fireshape}, or require the user to implement these algorithms themselves \cite{Gangl2020Fully}. Our approach leads to unique features, such as the possibility of discretizing and solving the state and adjoint systems differently as well as the choice of the scalar product for the computation of the (shape) gradients, and also gives rise to mesh independent behavior, as shown in Section~\ref{sec:examples}. Moreover, cashocs is the only one of these packages that has implemented a remeshing feature for shape optimization problems.

\subsection{Mathematical Background}
\label{ssec:background}

Let us begin with stating the general form of the optimization problems our software can solve. Optimal control problems have the form
\begin{linenomath*}
	\begin{equation}
		\label{eq:abstract_ocp}
		\min_{y, u}\ \costfunction(y, u) \quad \text{ s.t. } \quad e(y, u) = 0 \quad \text{ and } \quad u\in U_\text{ad},
	\end{equation}
\end{linenomath*}
where $u \in U$ and $y \in Y$ are the control and state variables, $U$ and $Y$ are appropriate Banach spaces, and the set of admissible controls $U_\text{ad} \subseteq U$ is used to model additional constraints on the control variable. Moreover, $\costfunction\colon Y \times U \to \R$ is the cost functional and $e\colon Y \times U \to Z^*$ is a PDE constraint, which we interpret in the following weak sense
\begin{linenomath*}
	\begin{equation*}
		\text{Find } y \in Y \text{ such that } \qquad \dual{e(y, u)}{p}{Z^*, Z} = 0 \qquad \text{ for all } p\in Z.
	\end{equation*}
\end{linenomath*}
Here, $Z^*$ denotes the topological dual space of $Z$, and $\dual{\varphi}{x}{Z^*, Z}$ denotes the duality pairing of $\varphi \in Z^*$ and $x\in Z$. 

Shape optimization problems have the form
\begin{linenomath*}
	\begin{equation}
		\label{eq:abstract_sop}
		\min_{y, \Omega}\ \costfunction(y, \Omega) \quad \text{ s.t. } \quad e(y, \Omega) = 0 \quad \text{ and } \quad \Omega \in \admissiblegeom,
	\end{equation}
\end{linenomath*}
where $y$ is again the state variable, and the set of admissible domains $\admissiblegeom$ is used to incorporate additional geometrical constraints. We interpret the PDE constraint $e(y,\Omega) = 0$ in the following weak sense
\begin{linenomath*}
	\begin{equation*}
		\text{Find } y \in Y(\Omega) \text{ such that } \quad \dual{e(y, \Omega)}{p}{Z(\Omega)^*, Z(\Omega)} = 0 \quad \text{ for all } p\in Z(\Omega).
	\end{equation*}
\end{linenomath*}
In particular, this means that the PDE constraint is given on the domain $\Omega$, and it is the shape of this domain that is subjected to optimization.

Problems~\eqref{eq:abstract_ocp} and~\eqref{eq:abstract_sop} are prototypes for the kinds of problems that cashocs can solve, and we refer the reader to Section~\ref{sec:examples} for illustrative examples. As mentioned previously, these kinds of problems are usually solved with the adjoint approach, whose derivation is beyond the scope of this paper. Hence, we refer the reader to \cite{Hinze2009Optimization, Troeltzsch2010Optimal} and \cite{Delfour2011Shapes, Schulz2016Efficient, Blauth2020Nonlinear} for a discussion and derivation of the adjoint approach for optimal control and shape optimization problems, respectively.

\section{Software description}
\label{sec:description}


\subsection{Software Architecture}
\label{ssec:architecture}

\begin{figure}[!t]
	\tikzstyle{block} = [rectangle, draw, text width=12em, text centered, rounded corners, minimum height=4em]
	\tikzstyle{splitblock} = [rectangle split, draw, text centered, rounded corners, minimum height=4em, rectangle split parts=2]
	\tikzstyle{line} = [draw, -latex']
	
	\begin{tikzpicture}[scale=0.2]
	\node[splitblock](config) at (0,0) {\textbf{configuration file}
		\nodepart[align=left]{two} optimization parameters \\ PDE parameters \\ output parameters};
	
	\node[splitblock, below of=config, node distance=10em](script) {\textbf{user script}
		\nodepart[align=left]{two} PDE constraint(s) \\ cost functional};
	
	\node[splitblock, right of=script, node distance=18em](optimization-problem) {\textbf{optimization problem}
		\nodepart[align=left]{two} \texttt{cashocs.OptimalControlProblem} \\ \texttt{cashocs.ShapeOptimizationProblem} };
	
	\node[splitblock](output) at (config -| optimization-problem) {\textbf{output}
		\nodepart[align=left]{two} (numerical) solution \\ history of the optimization \\ visualization };
	
	\path[line] (config) -- node[text width=2.5cm, midway, right=-0.5cm, align=center] {is \\loaded \\ into} (script);
	\path[line] (script) -- node[text width=2.5cm, midway, above, align=center] {calls} (optimization-problem);
	\path[line] (optimization-problem) -- node[text width=2.5cm, midway, right=-0.3cm, align=center] {generates} (output);
	\end{tikzpicture}
	\caption{Architecture of cashocs.}
	\label{fig:architecture}
\end{figure}
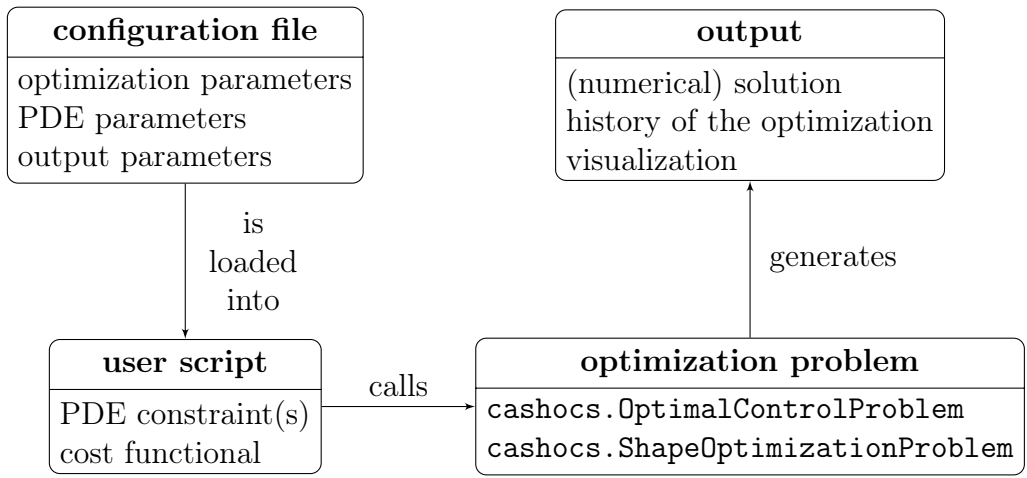

To solve optimization problems with cashocs, the user has to do the following. First, they have to implement the problem in a Python script, including the definition of the computational mesh, the state system, and the cost functional. To do so, they can use the same syntax as for defining the problem in FEniCS \cite{Logg2012Automated, Alnes2015FEniCS}, with only very minor modifications, resulting in a simple, high-level user interface that supports many important types of optimization problems. Second, the user has to define a configuration file that specifies the parameters for the solution of the state system and the optimization algorithm, which is loaded into the user script. Then, one can set up an optimization problem using \texttt{cashocs.OptimalControlProblem} or \texttt{cashocs.ShapeOptimizationProblem}, respectively, and solve it with the \texttt{solve} method of the respective class. Internally, our software utilizes the symbolic automatic differentiation capabilities of the Unified Form Language \cite{Logg2012Automated, Ham2019Automated} to compute the required (shape) derivatives and the variational formulation of the adjoint systems. Moreover, cashocs uses FEniCS to generate and compile \texttt{C++} code for the finite element assembly of the problems and PETSc \cite{Balay2020PETSc} is used to solve the arising linear systems, which makes the solution of the problems very efficient. A schematic overview of cashocs' architecture can be seen in Figure~\ref{fig:architecture}.


\subsection{Software Functionalities}
\label{ssec:functionalities}

Our software can treat linear and nonlinear systems of PDE constraints for steady state and transient conditions, as long as they can be implemented as (sequence of) variational formulation(s) in FEniCS. Further, cashocs deals with additional control constraints using projection techniques and can be used to solve state constrained problems, e.g., by means of a Moreau-Yosida regularization (see, e.g., \cite{Hinze2009Optimization}). We present two model problems constrained by Poisson's equation which illustrate the simplicity of cashocs' interface in Section~\ref{sec:examples} and refer the reader to the tutorial at \url{https://cashocs.readthedocs.io/en/latest/tutorial_index.html} for a detailed description of our software's capabilities for these more complex settings.

The following algorithms are available for shape optimization and optimal control problems in cashocs
\begin{itemize}
	\item the gradient descent method,
	
	\item nonlinear conjugate gradient methods (NCG) methods,
	
	\item limited memory BFGS (L-BFGS) methods.
\end{itemize}
Note, that particularly for shape optimization these algorithms correspond to the state of the art, with the L-BFGS methods being introduced in \cite{Schulz2016Efficient}, and the NCG methods in \cite{Blauth2020Nonlinear}. Additionally, the following optimization algorithms are available for optimal control problems only
\begin{itemize}
	\item a truncated Newton method,
	
	\item a primal-dual active set method.
\end{itemize}
Note, that for optimal control problems, all methods can also treat box constraints for the control variable using projection techniques. The user can adjust the behavior of these algorithms using the configuration file, where, e.g., the relative and absolute stopping tolerances, maximum number of iterations, and other, algorithm specific, parameters can be modified.

Additional features of cashocs include, among others, the possibility to use different discretizations for state and adjoint systems, the implementation of a Picard iteration for solving coupled systems, the possibility to specify which scalar product is used for the computation of the (shape) gradient, and remeshing for shape optimization problems, which utilizes the mesh generation software Gmsh \cite{Geuzaine2009Gmsh}.



\section{Illustrative Examples}
\label{sec:examples}

To demonstrate our software's simplicity for defining PDE constrained optimization problems as well as its efficiency for solving them, we now investigate two model problems, one for optimal control and one for shape optimization. Note, that a variety of other examples for using cashocs can be found in the tutorial at \url{https://cashocs.readthedocs.io/en/latest/tutorial_index.html}.

\subsection{Optimal Control}
\label{ssec:ex_optimal_control}

As a model optimal control problem we consider the following one from~\cite{Hinze2009Optimization}
\begin{linenomath*}
	\begin{equation}
		\label{eq:example_ocp}
		\begin{aligned}
			&\min\; \costfunction(y,u) = \frac{1}{2} \int_{\Omega} \left( y - y_d \right)^2 \text{ d}x + \frac{\alpha}{2} \int_{\Omega} u^2 \text{ d}x \\
			&\text{ subject to } \quad \left\lbrace \quad
			\begin{alignedat}{2}
			-\Delta y &= u \quad &&\text{ in } \Omega,\\
			y &= 0 \quad &&\text{ on } \Gamma = \partial\Omega.
			\end{alignedat} \right.
		\end{aligned}
	\end{equation}
\end{linenomath*}
This optimal control problem has a tracking-type cost functional with a Tikhonov regularization for the control variable. The PDE constraint is given by a Poisson equation with homogeneous Dirichlet boundary conditions, and the control variable $u$ enters the PDE as a right-hand side. The weak formulation of this PDE constraint is given by
\begin{linenomath*}
	\begin{equation}
		\label{eq:weak_poisson}
		\text{Find } y \in H^1_0(\Omega) \text{ such that } \quad \integral{\Omega} \grad y \cdot \grad p - u p \dmeas{x} = 0 \quad \text{ for all } p \in H^1_0(\Omega).
	\end{equation}
\end{linenomath*}
For this example, let us use $\Omega = (0,1)^2$, $\alpha = \num{1e-4}$, and
\begin{linenomath*}
	\begin{equation*}
		y_d(x) = x_1^2(1-x_1)\ x_2^2(1-x_2).
	\end{equation*}
\end{linenomath*}
For the discretization of the domain we use a uniform triangular mesh which divides $\Omega$ into $n\times n$ squares that are halved to create triangles. To solve this problem with cashocs, we can use the code shown in Listing~\ref{list_ocp}, which we briefly discuss in the following.

\begin{nolinenumbers}
\begin{python}[caption={Code for solving problem~\eqref{eq:example_ocp} with cashocs.}, captionpos=b, label={list_ocp}]
from fenics import *
import cashocs

# define mesh and volume measure
n = 64
mesh = UnitSquareMesh(n, n)
dx = Measure('dx', mesh)

# function space for linear Lagrange elements
V = FunctionSpace(mesh, 'CG', 1)
# define state, adjoint, and control variables
y = Function(V)
p = Function(V)
u = Function(V)

# define the weak form of the PDE constraint
e = inner(grad(y), grad(p))*dx - u*p*dx
# define the boundary conditions
bdry = CompiledSubDomain('on_boundary')
bcs = DirichletBC(V, Constant(0), bdry)

# define desired state and cost functional
x = SpatialCoordinate(mesh)
y_d = pow(x[0],2)*(1 - x[0])*pow(x[1],2)*(1-x[1])
J = 0.5*pow(y - y_d,2)*dx + 1e-4/2*pow(u,2)*dx

# solve the optimization problem with cashocs
cfg = cashocs.create_config('config.ini')
ocp = cashocs.OptimalControlProblem(e, bcs, J, y, u, p, cfg)
ocp.solve()
\end{python}
\end{nolinenumbers}
\begin{python}[float=t, caption={Minimal configuration file \texttt{config.ini} for problem~\eqref{eq:example_ocp}.}, captionpos=b, label={list_ocp_config}]
[OptimizationRoutine]
algorithm = ncg
rtol = 1e-3
maximum_iterations = 50

# additional parameters
# ...
\end{python}

Note, that as our software is based on FEniCS, we refer the reader to \cite[Chapter~1]{Logg2012Automated}, where the syntax of FEniCS is explained using several descriptive examples. In Listing~\ref{list_ocp}, we begin by importing FEniCS and cashocs in lines~1 and~2. Next, we define the mesh with the \texttt{UnitSquareMesh} function, and set up the volume measure for integration, in lines~5--7. Subsequently, we define a function space of linear Lagrange elements in line~10, and define the functions \texttt{y}, \texttt{p}, and \texttt{u}. These are used to define the weak form of the PDE constraint in line~17, where the function \texttt{p} plays the role of the test function.
\begin{figure}[!b]
	\centering
	\begin{subfigure}{0.4\textwidth}
		\includegraphics[width=\textwidth, trim=25cm 0cm 25cm 0cm, clip]{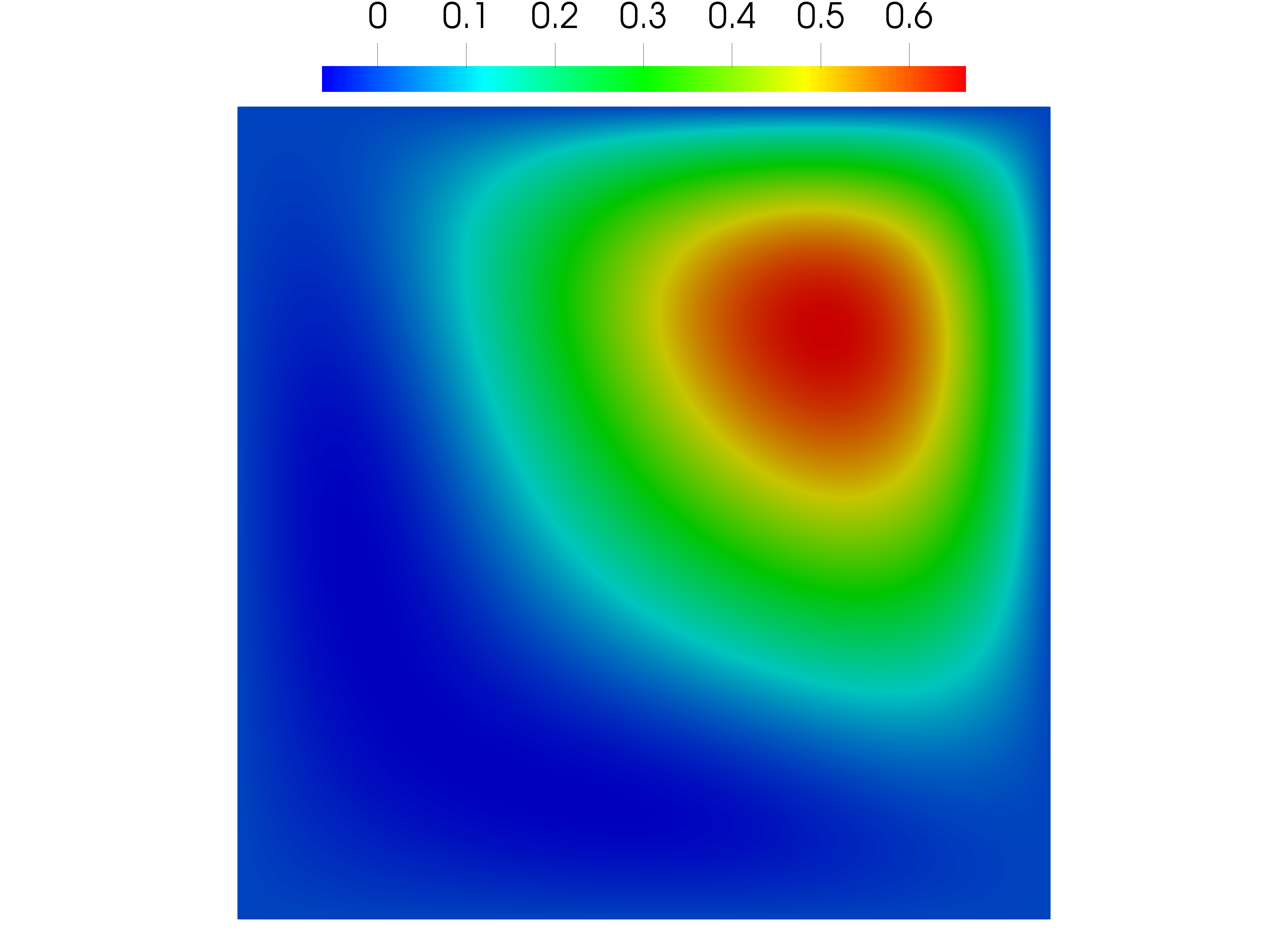}
		\caption{Optimal control $u$.}
	\end{subfigure}
	\hfil
	\begin{subfigure}{0.4\textwidth}
		\includegraphics[width=\textwidth, trim=25cm 0cm 25cm 0cm, clip]{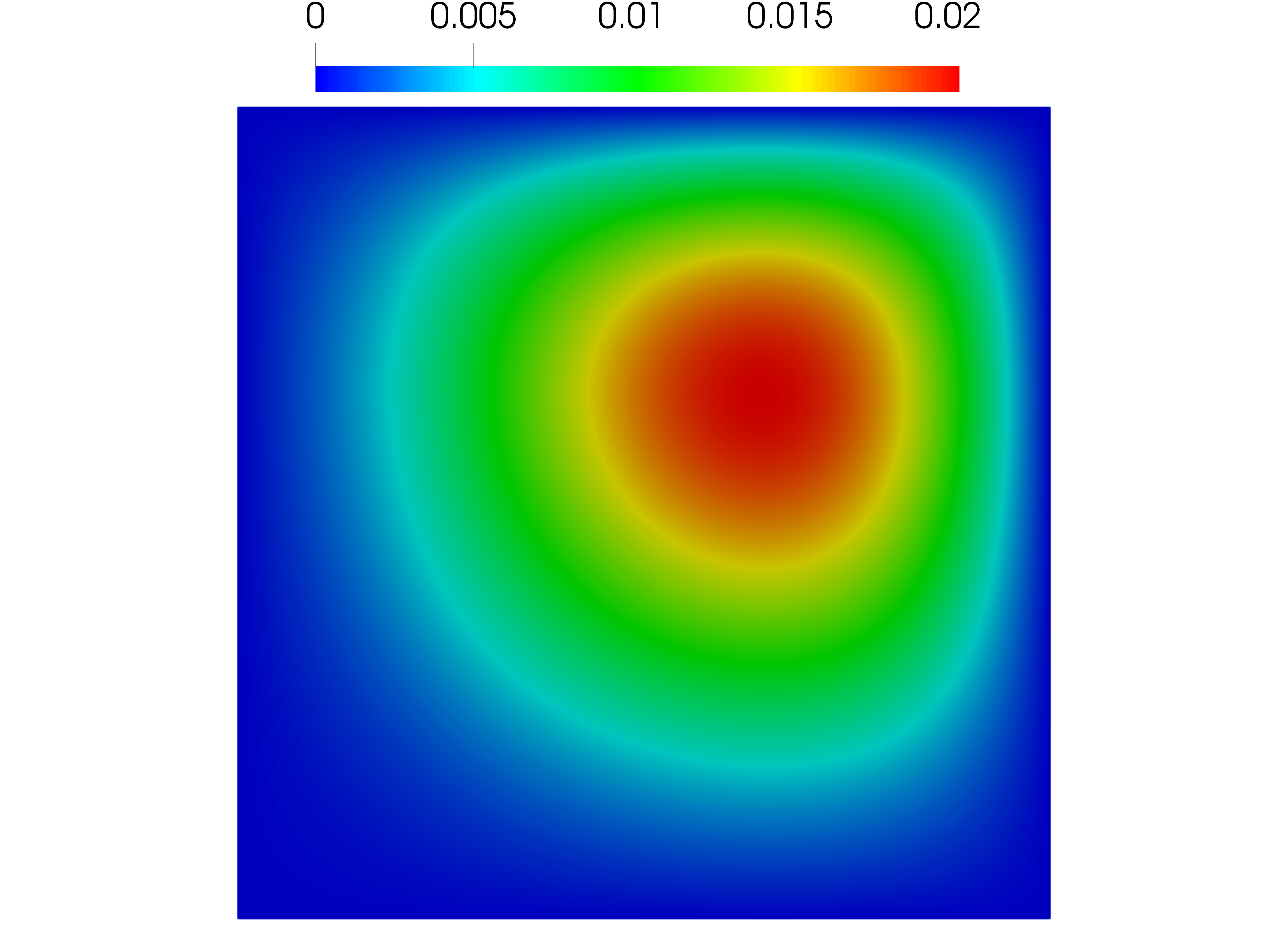}
		\caption{Optimal state $y$.}
	\end{subfigure}
	\caption{Numerical solution of problem~\eqref{eq:example_ocp}, computed with the Dai-Yuan NCG method.}
	\label{fig:sol_ocp}
\end{figure}
In lines~19 and~20 the Dirichlet boundary conditions for the Poisson problem are defined. Note, that up until now, we only used commands from FEniCS with the following minor variations. Instead of defining \texttt{y} as a \texttt{TrialFunction} and \texttt{p} as a \texttt{TestFunction}, both are now \texttt{Function} objects. Additionally, instead of defining the (linear) PDE constraint using its left- and right-hand sides, we define it as we would for a nonlinear variational problem in FEniCS, analogously to the form in \eqref{eq:abstract_ocp} and \eqref{eq:weak_poisson}. In lines~23 and~24 we define the desired state, which is used in line~25 to define the cost functional. Again, we have only used FEniCS commands for these operations. To invoke cashocs to solve this problem, all we have to do is loading the configuration file into the script in line~28, initializing the \texttt{OptimalControlProblem} in line~29, and calling its solve method subsequently. In total, we have to add only three additional lines of code to solve the problem. Note, that a minimal configuration file for the code is shown in Listing~\ref{list_ocp_config}, and for a detailed description of the configuration files for optimal control problems we refer to \url{https://cashocs.readthedocs.io/en/latest/demos/optimal_control/doc_config.html}. Note, that the scalar product used for computing the gradient of the cost functional can be determined by the user, as is explained in the tutorial. The default configuration uses the $L^2(\Omega)$ scalar product which is suitable for our model problem.

A plot of the computed optimal control and state using the Dai-Yuan nonlinear CG method is shown in Figure~\ref{fig:sol_ocp}. Moreover, Table~\ref{tab:optimal_control_mesh_indep} shows the amount of iterations the optimization algorithms need to solve this problem for a sequence of finer meshes, using $n=$ 16, 32, 64, and 128 subdivisions. We observe that all algorithms show mesh independent behavior as they basically need the same number of iterations for convergence regardless of the discretization.

\begin{table}[!t]
	\centering
	{\footnotesize
		\begin{tabular}{r r r r r r}
			\toprule
			$n$ & \hspace{0.5em} & {GD} & {NCG} & {L-BFGS} & {Newton} \\
			\midrule
			16 & & 32 & 10 & 6 & 1 \\
			\midrule
			32 & & 33 & 10 & 6 & 1 \\
			\midrule
			64 & & 33 & 10 & 6 & 1 \\
			\midrule
			128 & & 33 & 10 & 6 & 1 \\
			\bottomrule
		\end{tabular}
		\caption{Required iterations to reach the stopping criterion for problem~\eqref{eq:example_ocp}.}
		\label{tab:optimal_control_mesh_indep}
	}
\end{table}

\subsection{Shape Optimization}
\label{ssec:ex_shape_optimization}

As model problem for shape optimization we consider the following one from \cite{Blauth2020Nonlinear, Etling2020First}
\begin{linenomath*}
	\begin{equation}
		\label{eq:example_shape}
		\begin{aligned}
			&\min_{y, \Omega}\; \costfunction(y, \Omega) = \int_\Omega y \text{ d}x \\
			&\text{subject to} \quad \left\lbrace \quad
			\begin{alignedat}{2}
			-\Delta y &= f \quad &&\text{ in } \Omega,\\
			y &= 0 \quad &&\text{ on } \Gamma = \partial\Omega.
			\end{alignedat} \right.
		\end{aligned}
	\end{equation}
\end{linenomath*}
For this problem, the PDE constraint is, again, given by a Poisson problem with homogeneous Dirichlet boundary conditions, so that its weak form is given by \eqref{eq:weak_poisson} with $u$ replaced by $f$. 

We proceed analogously to \cite{Blauth2020Nonlinear, Etling2020First} and use as initial guess for the domain $\Omega_0$ the unit circle in $\mathbb{R}^2$, and for the right-hand side $f$ we use
\begin{linenomath*}
	\begin{equation*}
		f(x) = 2.5 \left( x_1 + 0.4 - x_2^2 \right)^2 + x_1^2 + x_2^2 - 1.
	\end{equation*}
\end{linenomath*}
We discretize $\Omega_0$ with a uniform triangular mesh by dividing the circle into $n$ smaller strips, which are then meshed uniformly. 
This problem can be solved with cashocs using the code provided in Listing~\ref{list_sop}, which we briefly discuss in the following. As before, we refer to \cite[Chapter~1]{Logg2012Automated} for a detailed introduction to the syntax of FEniCS, which we also use for the problem definition in cashocs.

\begin{nolinenumbers}
\begin{python}[caption={Code for solving problem~\eqref{eq:example_shape} with cashocs.}, captionpos=b, label={list_sop}]
from fenics import *
import cashocs

# define mesh and volume measure
n = 64
mesh = UnitDiscMesh.create(MPI.comm_world, n, 1, 2)
dx = Measure('dx', mesh)

# function space of linear Lagrange elements
V = FunctionSpace(mesh, 'CG', 1)
# state and adjoint variables
y = Function(V)
p = Function(V)

# right-hand side
x = SpatialCoordinate(mesh)
f = 2.5*pow(x[0] + 0.4 - pow(x[1], 2), 2) \
	+ pow(x[0], 2) + pow(x[1], 2) - 1
# define the PDE constraint
e = inner(grad(y), grad(p))*dx - f*p*dx
# define the boundary conditions
bdry = CompiledSubDomain('on_boundary')
mf_bdry = MeshFunction('size_t', mesh, dim=1)
bdry.mark(mf_bdry, 1)
bcs = DirichletBC(V, Constant(0), mf_bdry, 1)

# cost functional
J = y*dx

# solve the problem with cashocs
cfg = cashocs.create_config('config.ini')
sop = cashocs.ShapeOptimizationProblem(e, bcs, J, y, p, 
										mf_bdry, cfg)
sop.solve()
\end{python}
\end{nolinenumbers}

\begin{python}[float=!b, caption={Minimal configuration file \texttt{config.ini} for problem~\eqref{eq:example_shape}.}, captionpos=b, label={list_sop_config}]
[OptimizationRoutine]
algorithm = ncg
rtol = 5e-3
maximum_iterations = 50

[ShapeGradient]
shape_bdry_def = [1]
shape_bdry_fix = []

# additional parameters
# ...
\end{python}
The code is very similar to the one in Listing~\ref{list_ocp} as we again have a Poisson equation as PDE constraint. We start the script by importing FEniCS and cashocs. Then, we define the mesh and volume measure, now using the function \texttt{UnitDiscMesh}, in lines~5--7. For the discretization of the Poisson equation, we again use linear Lagrange elements whose corresponding function space is defined in line~10, and the functions \texttt{y} and \texttt{p} are defined in lines~12 and~13. Thereafter, we define the right-hand side of the Poisson problem, using \texttt{SpatialCoordinate} in lines~16--18, which is then used to define the weak form of the Poisson equation in line~20. As for optimal control problems, the only major differences to traditional FEniCS syntax are that \texttt{y} and \texttt{p} are \texttt{Function} objects, and that the PDE constraint is written in the sense of \eqref{eq:abstract_sop} and \eqref{eq:weak_poisson}. Subsequently, we set up a FEniCS \texttt{MeshFunction} for the boundaries, which is used to define the Dirichlet boundary conditions. Moreover, this is used to define which boundaries are fixed via the configuration file (cf.\ lines~7--8 of Listing~\ref{list_sop_config}). Finally, we define the cost functional in line~28. For solving this problem with cashocs, we proceed analogously to Listing~\ref{list_ocp}, and first load the configuration file, then set up the \texttt{ShapeOptimizationProblem}, and finally call its \texttt{solve} method in lines~31--34. Note, that a minimal configuration file for this problem is shown in Listing~\ref{list_sop_config}. A detailed discussion of the configuration files for shape optimization can be found at \url{https://cashocs.readthedocs.io/en/latest/demos/shape_optimization/doc_config.html}.

Note, that the scalar product used for computing the shape gradient is based on the linear elasticity equations (see, e.g., \cite{Schulz2016Computational, Blauth2020Nonlinear, Etling2020First}). The corresponding bilinear form $a$ is given by
\begin{equation*}
	a(\vectorfield, \mathcal{W}) = \integral{\Omega} 2 \mu\ \varepsilon(\vectorfield) : \varepsilon(\mathcal{W}) + \lambda\ \text{div}(\vectorfield)\ \text{div}(\mathcal{W}) + \delta\ \vectorfield \cdot \mathcal{W} \dmeas{x},
\end{equation*}
where $\varepsilon(\vectorfield) = \nicefrac{1}{2} (D\vectorfield + D\vectorfield\transposed)$ is the symmetric part of the Jacobian. The default values for the elasticity parameters are $\mu = 1$, $\lambda = 0$, and $\delta = 0$ and can be altered through the configuration file.

\begin{figure}[!t]
	\centering
	\includegraphics[width=0.4\textwidth]{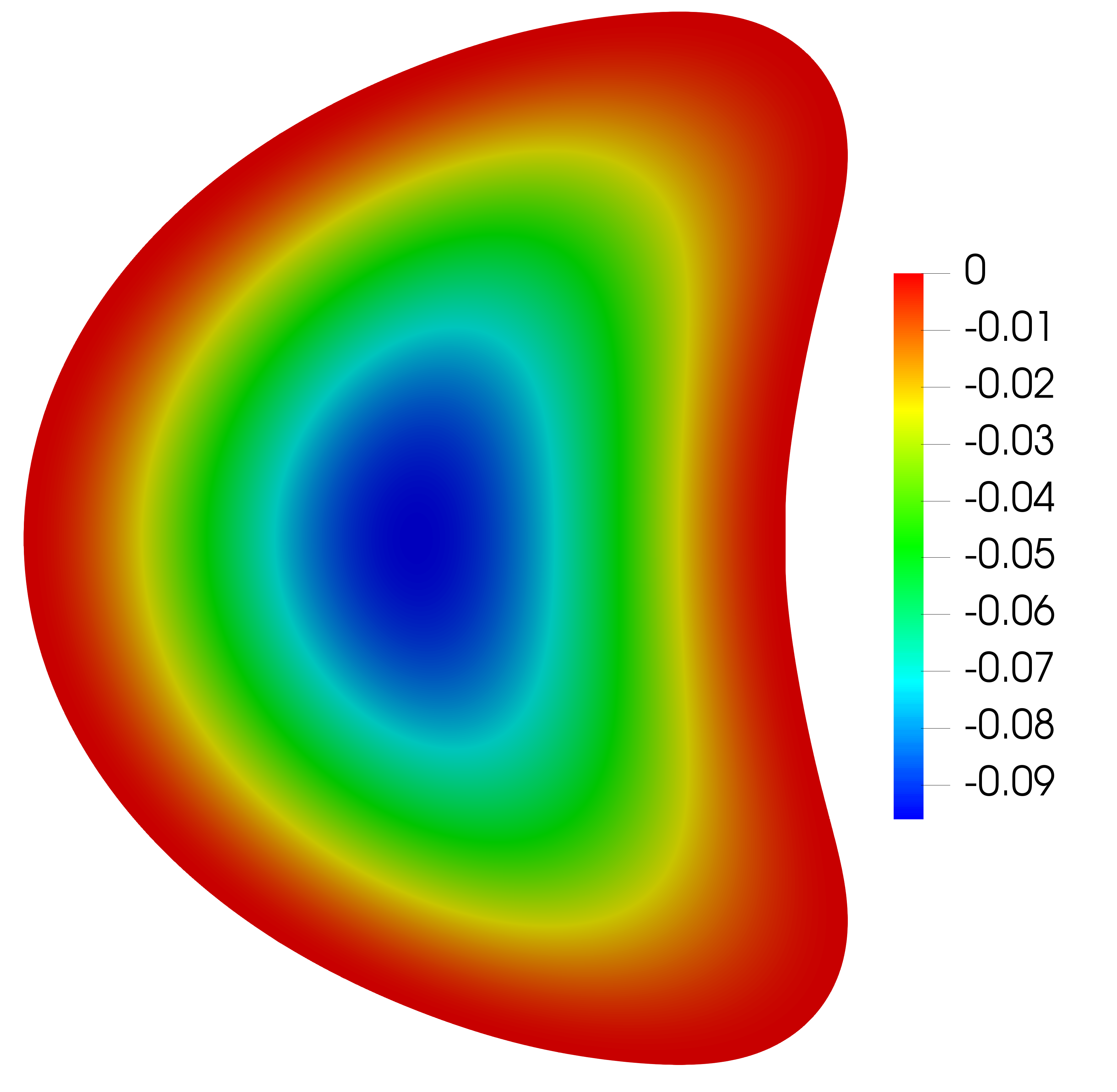}
	\caption{Numerical solution of problem~\eqref{eq:example_shape}, optimal state $y$ on the optimal domain $\Omega$, computed with the Dai-Yuan NCG method.}
	\label{fig:sol_sop}
\end{figure}

\begin{table}[!b]
	\centering
	{\footnotesize
		\begin{tabular}{r r r r r}
			\toprule
			$n$ & \hspace{0.5em} & {GD} & {NCG} & {L-BFGS} \\
			\midrule
			16 & & 46 & 20 & 12 \\
			\midrule
			32 & & 47 & 19 & 11 \\
			\midrule
			64 & & 47 & 19 & 11 \\
			\midrule
			128 & & 47 & 19 & 11 \\
			\bottomrule
		\end{tabular}
		\caption{Required iterations to reach the stopping criterion for problem~\eqref{eq:example_shape}.}
		\label{tab:shape_optimization_mesh_indep}
	}
\end{table}

A plot of the optimal state on the optimal domain, computed with the Dai-Yuan NCG method in cashocs, is given in Figure~\ref{fig:sol_sop}. Moreover, we also show the number of iterations required by the algorithms on successively finer discretizations of $n=$ 16, 32, 64, and 128 strips for the unit circle in Table~\ref{tab:shape_optimization_mesh_indep}. As before, we see that the number of iterations basically stays the same regardless of the discretization, which shows that we also have mesh independent behavior for shape optimization problems.

\section{Impact}
\label{sec:impact}

Our software enables users to treat complex, coupled, and highly nonlinear PDE constrained optimization problems in an automated fashion. The user is only required to define the PDE constraint and cost functional using basically the same syntax as for defining these objects in FEniCS. Thanks to the high-level user interface, the corresponding optimization problem can then be solved by adding only three additional lines of code. 
Our approach of implementing a discretization of the continuous adjoint approach leads to mesh independent behavior of the optimization algorithms, as shown in Section~\ref{sec:examples}, making our software attractive for science and industry. In fact, cashocs has already been used to treat highly nonlinear optimization problems for parameter identification and optimal control in the context of chemical microreactors in \cite{Blauth2020Optimal}. It has also been used in \cite{Blauth2020Nonlinear} for a numerical benchmark of NCG methods for shape optimization. Moreover, cashocs is used at Fraunhofer ITWM to solve PDE constrained optimization problems for industrial applications. Due to the generality of our software, which can treat lots of important classes of cost functionals and PDE constraints, it can be applied to many relevant problems in science and industry, automating their solution in an efficient and user friendly way.



%
%
%

\section{Conclusions}
\label{sec:conclusions}

We have presented cashocs, a software for numerically solving PDE constrained shape optimization and optimal control problems. The software automatically derives the required adjoint systems and (shape) derivatives, and implements a discretization of the continuous adjoint approach. Our software inherits FEniCS' high-level user interface which allows for a straightforward definition and solution of PDE constrained optimization problems. Additionally, the user still retains control over many important parameters for the optimization, ranging from the solution of the PDEs to the optimization algorithm, which allows them to make precise adjustments to the numerical solution of their problems. 


\section{Conflict of Interest}

We wish to confirm that there are no known conflicts of interest associated with this publication and there has been no significant financial support for this work that could have influenced its outcome.

\section*{Acknowledgements}
\label{ack}

The author gratefully acknowledges financial support from the Fraunhofer Institute for Industrial Mathematics ITWM.

\bibliographystyle{elsarticle-num} 
\bibliography{lit.bib}

\begin{thebibliography}{10}
\expandafter\ifx\csname url\endcsname\relax
  \def\url#1{\texttt{#1}}\fi
\expandafter\ifx\csname urlprefix\endcsname\relax\def\urlprefix{URL }\fi
\expandafter\ifx\csname href\endcsname\relax
  \def\href#1#2{#2} \def\path#1{#1}\fi

\bibitem{Blauth2020Optimal}
S.~Blauth, C.~Leith\"{a}user, R.~Pinnau, {O}ptimal {C}ontrol of the {S}abatier
  {P}rocess in {M}icrochannel {R}eactors (2020).
\newblock \href {http://arxiv.org/abs/2007.12457} {\path{arXiv:2007.12457}}.

\bibitem{Pinnau2004Optimal}
R.~Pinnau, G.~Th\"{o}mmes, Optimal boundary control of glass cooling processes,
  Math. Methods Appl. Sci. 27~(11) (2004) 1261--1281.
\newblock \href {http://dx.doi.org/10.1002/mma.500}
  {\path{doi:10.1002/mma.500}}.

\bibitem{Hinze2002optimal}
M.~Hinze, R.~Pinnau, An optimal control approach to semiconductor design, Math.
  Models Methods Appl. Sci. 12~(1) (2002) 89--107.
\newblock \href {http://dx.doi.org/10.1142/S0218202502001568}
  {\path{doi:10.1142/S0218202502001568}}.

\bibitem{Blauth2020Shape}
S.~Blauth, C.~Leith\"{a}user, R.~Pinnau, Shape sensitivity analysis for a
  microchannel cooling system, J. Math. Anal. Appl. 492~(2) (2020) 124476.
\newblock \href {http://dx.doi.org/10.1016/j.jmaa.2020.124476}
  {\path{doi:10.1016/j.jmaa.2020.124476}}.

\bibitem{Schmidt2013Three}
S.~Schmidt, C.~Ilic, V.~Schulz, N.~R. Gauger, {T}hree-{D}imensional
  {L}arge-{S}cale {A}erodynamic {S}hape {O}ptimization {B}ased on {S}hape
  {C}alculus, AIAA Journal 51~(11) (2013) 2615--2627.
\newblock \href {http://dx.doi.org/10.2514/1.J052245}
  {\path{doi:10.2514/1.J052245}}.

\bibitem{Gangl2015Shape}
P.~Gangl, U.~Langer, A.~Laurain, H.~Meftahi, K.~Sturm, Shape {O}ptimization of
  an {E}lectric {M}otor {S}ubject to {N}onlinear {M}agnetostatics, SIAM J. Sci.
  Comput. 37~(6) (2015) B1002--B1025.
\newblock \href {http://dx.doi.org/10.1137/15100477X}
  {\path{doi:10.1137/15100477X}}.

\bibitem{Mitusch2019dolfin}
S.~K. Mitusch, S.~W. Funke, J.~S. Dokken, dolfin-adjoint 2018.1: automated
  adjoints for {FEniCS} and {F}iredrake, Journal of Open Source Software 4~(38)
  (2019) 1292.
\newblock \href {http://dx.doi.org/10.21105/joss.01292}
  {\path{doi:10.21105/joss.01292}}.

\bibitem{Dokken2020Automatic}
J.~S. Dokken, S.~K. Mitusch, S.~W. Funke, Automatic shape derivatives for
  transient {PDEs} in {FEniCS} and {Firedrake} (2020).
\newblock \href {http://arxiv.org/abs/2001.10058} {\path{arXiv:2001.10058}}.

\bibitem{Paganini2020Fireshape}
A.~Paganini, F.~Wechsung, Fireshape: a shape optimization toolbox for
  {F}iredrake (2020).
\newblock \href {http://arxiv.org/abs/2005.07264} {\path{arXiv:2005.07264}}.

\bibitem{Gangl2020Fully}
P.~Gangl, K.~Sturm, M.~Neunteufel, J.~Schöberl, Fully and semi-automated shape
  differentiation in {NGSolve}, Struct. Multidiscip. Optim.\href
  {http://dx.doi.org/10.1007/s00158-020-02742-w}
  {\path{doi:10.1007/s00158-020-02742-w}}.

\bibitem{Hinze2009Optimization}
M.~Hinze, R.~Pinnau, M.~Ulbrich, S.~Ulbrich, Optimization with {PDE}
  constraints, Vol.~23 of Mathematical Modelling: Theory and Applications,
  Springer, New York, 2009.
\newblock \href {http://dx.doi.org/10.1007/978-1-4020-8839-1}
  {\path{doi:10.1007/978-1-4020-8839-1}}.

\bibitem{Troeltzsch2010Optimal}
F.~Tr\"{o}ltzsch, Optimal {C}ontrol of {P}artial {D}ifferential {E}quations,
  Vol. 112 of Graduate Studies in Mathematics, American Mathematical Society,
  Providence, RI, 2010.
\newblock \href {http://dx.doi.org/10.1090/gsm/112}
  {\path{doi:10.1090/gsm/112}}.

\bibitem{Delfour2011Shapes}
M.~C. Delfour, J.-P. Zol\'{e}sio, Shapes and {G}eometries, 2nd Edition, Vol.~22
  of Advances in Design and Control, Society for Industrial and Applied
  Mathematics (SIAM), Philadelphia, PA, 2011.
\newblock \href {http://dx.doi.org/10.1137/1.9780898719826}
  {\path{doi:10.1137/1.9780898719826}}.

\bibitem{Schulz2016Efficient}
V.~H. Schulz, M.~Siebenborn, K.~Welker, Efficient {PDE} {C}onstrained {S}hape
  {O}ptimization {B}ased on {S}teklov-{P}oincar\'{e}-{T}ype {M}etrics, SIAM J.
  Optim. 26~(4) (2016) 2800--2819.
\newblock \href {http://dx.doi.org/10.1137/15M1029369}
  {\path{doi:10.1137/15M1029369}}.

\bibitem{Blauth2020Nonlinear}
S.~Blauth, {N}onlinear {C}onjugate {G}radient {M}ethods for {PDE} {C}onstrained
  {S}hape {O}ptimization {B}ased on {S}teklov-{P}oincaré-{T}ype {M}etrics
  (2020).
\newblock \href {http://arxiv.org/abs/2007.12891} {\path{arXiv:2007.12891}}.

\bibitem{Logg2012Automated}
A.~Logg, K.-A. Mardal, G.~N. Wells, et~al., Automated Solution of Differential
  Equations by the Finite Element Method, Springer, 2012.
\newblock \href {http://dx.doi.org/10.1007/978-3-642-23099-8}
  {\path{doi:10.1007/978-3-642-23099-8}}.

\bibitem{Alnes2015FEniCS}
M.~S. Aln{\ae}s, J.~Blechta, J.~Hake, A.~Johansson, B.~Kehlet, A.~Logg,
  C.~Richardson, J.~Ring, M.~E. Rognes, G.~N. Wells, {T}he {FEniCS} {P}roject
  {V}ersion 1.5, Archive of Numerical Software 3~(100).
\newblock \href {http://dx.doi.org/10.11588/ans.2015.100.20553}
  {\path{doi:10.11588/ans.2015.100.20553}}.

\bibitem{Ham2019Automated}
D.~A. Ham, L.~Mitchell, A.~Paganini, F.~Wechsung, Automated shape
  differentiation in the {U}nified {F}orm {L}anguage, Struct. Multidiscip.
  Optim. 60~(5) (2019) 1813--1820.
\newblock \href {http://dx.doi.org/10.1007/s00158-019-02281-z}
  {\path{doi:10.1007/s00158-019-02281-z}}.

\bibitem{Balay2020PETSc}
S.~Balay, S.~Abhyankar, M.~F. Adams, J.~Brown, P.~Brune, K.~Buschelman,
  L.~Dalcin, A.~Dener, V.~Eijkhout, W.~D. Gropp, D.~Karpeyev, D.~Kaushik, M.~G.
  Knepley, D.~A. May, L.~C. McInnes, R.~T. Mills, T.~Munson, K.~Rupp, P.~Sanan,
  B.~F. Smith, S.~Zampini, H.~Zhang, H.~Zhang,
  \href{https://www.mcs.anl.gov/petsc}{{PETS}c users manual}, Tech. Rep.
  ANL-95/11 - Revision 3.13, Argonne National Laboratory (2020).
\newline\urlprefix\url{https://www.mcs.anl.gov/petsc}

\bibitem{Geuzaine2009Gmsh}
C.~Geuzaine, J.-F. Remacle, Gmsh: {A} 3-{D} finite element mesh generator with
  built-in pre- and post-processing facilities, Internat. J. Numer. Methods
  Engrg. 79~(11) (2009) 1309--1331.
\newblock \href {http://dx.doi.org/10.1002/nme.2579}
  {\path{doi:10.1002/nme.2579}}.

\bibitem{Etling2020First}
T.~Etling, R.~Herzog, E.~Loayza, G.~Wachsmuth, First and {S}econd {O}rder
  {S}hape {O}ptimization {B}ased on {R}estricted {M}esh {D}eformations, SIAM J.
  Sci. Comput. 42~(2) (2020) A1200--A1225.
\newblock \href {http://dx.doi.org/10.1137/19M1241465}
  {\path{doi:10.1137/19M1241465}}.

\bibitem{Schulz2016Computational}
V.~Schulz, M.~Siebenborn, {C}omputational {C}omparison of {S}urface {M}etrics
  for {PDE} {C}onstrained {S}hape {O}ptimization, Comput. Methods Appl. Math.
  16~(3) (2016) 485--496.
\newblock \href {http://dx.doi.org/10.1515/cmam-2016-0009}
  {\path{doi:10.1515/cmam-2016-0009}}.

\end{thebibliography}

\end{document}